\documentclass{amsart}
\usepackage{graphicx}
\usepackage{amssymb}
\usepackage{enumitem}

\newtheorem{thm}{Theorem}[section]

\theoremstyle{definition}

\theoremstyle{remark}
\newtheorem{rem}[thm]{Remark}
\newtheorem{quest}[thm]{Question}

\begin{document}

\title[large singular sets]{Solutions to the minimal surface system with large singular sets}
\author{Connor Mooney}
\address{Department of Mathematics, UC Irvine}
\email{\tt mooneycr@uci.edu}
\author{Ovidiu Savin}
\address{Department of Mathematics, Columbia University}
\email{\tt savin@math.columbia.edu}

\begin{abstract}
Lawson and Osserman proved that the Dirichlet problem for the minimal surface system is not always solvable in the class of Lipschitz maps. However, it is known that minimizing sequences (for area) of Lipschitz graphs converge to objects called Cartesian currents. Essentially nothing is known about these limits. We show that such limits can have surprisingly large interior vertical and non-minimal portions. This demonstrates a striking discrepancy between the parametric and non-parametric area minimization problems in higher codimension. Moreover, our construction has the smallest possible dimension ($n = 3$) and codimension $(m = 2)$.
\end{abstract}
\maketitle

\section{Introduction}
It is well-understood what happens when we minimize area among Lipschitz graphs in codimension one. Namely, fix a smooth bounded domain $\Omega \subset \mathbb{R}^n$ and a smooth real-valued function $\varphi$ on $\partial \Omega$. If $\Omega$ is mean-convex, then in the class $\mathcal{A}$ of Lipschitz functions with boundary data $\varphi$, there is a unique minimizer $u$ of the area functional. Moreover, $u \in C^{\infty}\left(\overline{\Omega}\right)$, and the graph of $u$ minimizes area among all (not necessarily graphical) competitors in $\Omega \times \mathbb{R}$ (and all competitors in $\mathbb{R}^{n+1}$ if $\Omega$ is convex). If $\Omega$ is not mean-convex, then there need not exist a minimizer of area in $\mathcal{A}$. Instead, minimizing sequences (for area) in $\mathcal{A}$ will converge in $L^1(\Omega)$ to a function $u \in C^{\infty}(\Omega) \cap BV(\Omega)$ such that the area of the union of the graph of $u$ and the part of the cylinder $\partial \Omega \times \mathbb{R}$ between the graphs of $u$ and $\varphi$ is the infimum of areas of functions in $\mathcal{A}$. That is, boundary discontinuities can appear in the limit, but interior discontinuities cannot. The canonical example is the case that $\Omega$ is an annulus, and $\varphi$ is zero on the outer sphere and a large constant on the inner sphere. In that case, the graph of $u$ is part of a catenoid that becomes vertical on the inner boundary. See e.g. Giusti's book \cite{Giu}, Part II, and the references therein.

In \cite{LO} Lawson and Osserman investigated the extent to which similar results might hold in higher codimension. Fix $m \geq 2$, a smooth bounded domain $\Omega \subset \mathbb{R}^n$, and a smooth map $\varphi : \partial \Omega \rightarrow \mathbb{R}^m$. They showed that if $\Omega$ is convex and $n = 2$, then there is a map $u \in C^{\infty}\left(\overline{\Omega}\right)$ with values in $\mathbb{R}^m$ and boundary data $\varphi$ whose graph minimizes area among competitors in $\mathbb{R}^{n+m}$. The idea is to show that any Douglas-Rad\'{o} solution to the Plateau problem is graphical. However, they also showed that this map is not always unique. In higher dimension the results were even more discouraging. Let $H$ be the Hopf map from $\mathbb{S}^3$ to $\mathbb{S}^2$, given by
$$H(z_1,\,z_2) := (|z_1|^2 - |z_2|^2,\, 2z_1\overline{z_2}),\, z_1,\,z_2 \in \mathbb{C}.$$
It was shown for $R > 0$ large that there is no Lipschitz map $u$ from $B_1 \subset \mathbb{R}^4$ to $\mathbb{R}^3$ whose graph is a critical point of area in $\mathbb{R}^7$, such that $u|_{\partial B_1} = RH$. Moreover, the graph of the one-homogeneous Lipschitz map
$$u(x) = \frac{\sqrt{5}}{2}|x|H(x/|x|)$$
is a critical point of area in $\mathbb{R}^7$ with a singularity at the origin. (It is in fact an area-minimizing cone, as shown in the work of Harvey-Lawson on calibrated geometries \cite{HL1}). Thus, the uniqueness, existence, and regularity results known in codimension one on convex domains fail spectacularly in higher codimension.

However, in the spirit of what was done in codimension one for non-convex domains, it is reasonable to ask:
\begin{quest}\label{Quest}
What is the nature of the limiting objects that arise when we minimize area among graphs of Lipschitz maps with fixed boundary data in higher codimension?
\end{quest} 
\noindent Question \ref{Quest} is not only pertinent for problems involving the area functional, but also for a variety of calculus of variations problems where minimizers do not exist in traditional spaces (more on this below).

As a starting point to study Question \ref{Quest}, it is instructive to start specific: Take $\Omega = B_1 \subset \mathbb{R}^4,\, m = 3,$ and $\varphi = RH$. In \cite{DY}, Ding-Yuan proved that for any $R > 0$, there a smooth solution to the minimal surface system in $B_1 \backslash \{0\}$ of the form
$$u(x) = f(|x|)H(x/|x|)$$
with boundary data $RH$. When $R < \sqrt{5}/2$, $f$ is an analytic even function such that $f(0) = 0$, and $u$ is analytic in $B_1$. When $R = \sqrt{5}/{2}$, $f(s) = \sqrt{5}s/2$, and the graph of $u$ is a minimal cone with an isolated singularity at the origin. Finally, when $R > \sqrt{5}/2$, $f$ is an analytic even function and $f(0) > 0$. In all cases, the graph of $u$ is area-minimizing. The last case $R > \sqrt{5}/2$ deserves some explanation, since in that case $u$ is discontinuous at the origin. One can view $u$ as an $L^1$ limit of Lipschitz maps $u_{\epsilon}$ obtained by cutting $f$ off near zero, e.g. by linear interpolation on $[0,\,\epsilon]$. It is not hard to check that the graph of $u_{\epsilon}$ over $B_{\epsilon}$ has area of order $\epsilon$, so the areas of the graphs of $u_{\epsilon}$ converge to the area of the graph of $u$ over $B_1 \backslash \{0\}$. It is useful to think of the graph of $u$ has having a ``vertical piece" of zero area in $\{0\} \times \mathbb{R}^3$, which is a sort of limit of the graphs of $u_{\epsilon}$ over $B_{\epsilon}$. This piece captures the nontrivial topology of the Hopf map (which ensures that any continuous map in $\overline{B_1}$ with boundary data $RH$ takes the value $0$ somewhere). However, the graph of $u$ has zero density at points on this vertical piece. If we discard it, we are left with an analytic area-minimizing submanifold of $\mathbb{R}^7$ (the closure of the graph of $u$ over $B_1 \backslash \{0\}$) that avoids the origin. Roughly, the topology changes at $R = \sqrt{5}/2$ to accommodate area minimization. 

On the other hand, one can study Question \ref{Quest} from a general viewpoint, and try to apply direct methods in the calculus of variations. This was done by Giaquinta-Modica-Sou\v{c}ek in \cite{GMS1}, \cite{GMS2}, \cite{GMS3}, \cite{GMS4}. One can view graphs of Lipschitz maps as rectifiable currents. Sequences of such graphs with uniformly bounded area and $L^1$ norm have weak subsequential limits called Cartesian currents, which are rectifiable currents with additional structure. Roughly, they are the sum of two parts: one given by the graph of a BV map, and another part being purely vertical. Cartesian currents are not only useful for studying variational problems involving area, but also for studying problems in elasticity where fractures (discontinuities, captured by purely vertical parts in Cartesian currents) occur. We will briefly review the theory of Cartesian currents in Section \ref{Prelims}.

Minimizing sequences (for area) of Lipschitz graphs with fixed boundary will thus converge to a Cartesian current. By the lower semicontinuity of area (see Section \ref{Prelims}), this current has smaller area than any Lipschitz graph with the same boundary. It is natural to ask whether it minimizes mass among all rectifiable currents with the same boundary, has generalized mean curvature equal to zero, and has vertical part with no mass (such questions were in fact asked at the end of Chapter 6 in the book \cite{GMS2} of Giaquinta-Modica-Sou\v{c}ek). For the example discussed above, the answers are all ``yes," and it is tempting to guess that such properties hold more generally. In this paper we show that the answer to each of these questions is in fact ``no". That is, surprisingly large interior vertical and non-minimal patches can arise in the limit.

We now explain our result more precisely. We will first construct a map 
$$w = (w^1,\,w^2)$$ 
from a bounded convex domain $\Omega \subset \mathbb{R}^3$ to $\mathbb{R}^2$ such that $\partial_1w^1 \geq 0$, and $\{\partial_1w^1 = 0\}$ is the closure of an analytic convex sub-domain $\Omega_0 \subset \subset \Omega$. The map $w$ is analytic in $\Omega \backslash \overline{\Omega_0}$ and $\Omega_0$, $C^{1,1}$ in $\Omega$, solves the minimal surface system in $\Omega \backslash \overline{\Omega_0}$, and does not solve the minimal surface system in $\Omega_0$. We then swap the $x_1$ and $w^1$ axes, and represent the graph of $w$ by that of a new map $u$ defined by
$$u^1(y(x)) = x_1,\, u^2(y(x)) = w^2(x), \quad y(x) := (w^1(x),\,x_2,\,x_3).$$
The set 
$$\Sigma := y\left(\overline{\Omega_0}\right)$$ 
is an analytic embedded surface with boundary that resembles a potato chip, and 
$$U := y(\Omega)$$ 
is a small neighborhood of $\Sigma$. The map $u$ is an analytic solution to the minimal surface system in $U \backslash \Sigma$, and is discontinuous on $\Sigma$. We understand the graph of $u$ as being a rotation in $\mathbb{R}^5$ of the graph of $w$, with a vertical part of dimension three (a rotation of the graph of $w$ over $\Omega_0$) that projects down to $\Sigma$ (in fact, each point in the interior of $\Sigma$ is the projection to $U$ of part of an analytic curve). We prove:

\begin{thm}\label{main}
The map $w$ described above can be constructed such that its graph has smaller area than the graph of any Lipschitz map from $U$ to $\mathbb{R}^2$ with the same boundary data as $u$ on $\partial U$. 
\end{thm}

It is easy to see that the graph of $u$ is a Cartesian current given by the limit of a minimizing (for area) sequence of graphs of Lipschitz maps, obtained by rotating the graphs of $w + (x_1/k,\,0)$ and gluing these to the graph of $u$ near $\partial U$. Since the graph of $u$ is not minimal and has an interior vertical part with positive area, we get negative answers to the questions asked above.

It is worth emphasizing that the vertical part of our example occurs away from the boundary, which is the only place where such parts could appear in codimension one. In codimension one, the limit in fact solves an obstacle problem, where the obstacle is the cylinder over the boundary. We construct $w$ by solving a different free boundary problem (see discussion below). Thus, in higher codimension, graphicality can be seen as a constraint that produces interior free boundaries when minimizing area. We also emphasize that our construction cannot be carried out in dimension $n = 2$, essentially because a maximum principle holds for the gradient of solutions to the minimal surface system in that case. See Remarks \ref{MP} and \ref{MP2}. Our example thus has the smallest possible dimension and codimension.

\begin{rem}
The space of Cartesian currents enjoys good compactness properties, which quickly give the existence of minimizers of area in the class of Cartesian currents with fixed boundary. Our proof in fact shows that the graph of $u$ is such a minimizer, see Remark \ref{CC}. Theorem \ref{main} can thus be viewed as saying that minimizers of area in the class of Cartesian currents can have vertical parts of maximal dimension, that moreover are not minimal, and project to hypersurfaces in the domain. It would be interesting to prove some general regularity results for such minimizers. For example, one could hope to say the the vertical parts project to a set of codimension one, that off of this set and a smaller one of codimension four the minimizer is a smooth graph, that the vertical patches have some kind of boundary regularity, and that the minimizer is regular as a geometric object away from a small set. As a reference point, our example is $C^{1,1}$ (and no better) as a geometric object.

As a starting point one could ask whether a monotonicity formula holds for minimizers of area among Cartesian currents. Our example is not minimal, showing that the traditional one cannot hold (the trouble in the proof being that the cone over the intersection of a graph and a sphere centered on the graph need not be graphical), but this doesn't rule out the validity of an alternative version.
\end{rem}

\begin{rem}\label{BigNHood}
We construct $w$ using Cauchy-Kovalevskaya. As a result, 
the vertical and non-minimal part of the graph of $u$ is very close to the boundary. It would be interesting to produce examples where this part is far from the boundary (``comes from infinity"), and to understand its stability. We touch on a possible approach below, which we leave for later work.
\end{rem}

To motivate our construction, we will first build an example with a point singularity. Namely, we find a solution $w = (w^1,\,w^2)$ to the minimal surface system near the origin in $\mathbb{R}^3$ such that $\partial_1w^1 \geq 0$ and $\{\partial_1w^1 = 0\} = \{0\}$, and then perform the rotation described above. Similar examples were previously constructed in dimension and codimension $n = m = 3$, and moreover are special Lagrangian, see the papers of Nadirashvili-Vladut \cite{NV2} and Wang-Yuan \cite{WY2}. To our knowledge, ours is the first example of a singular solution to the minimal surface system in the smallest possible dimension and codimension. 

To get examples with larger singularities, the idea is to rotate the graph of the one-point example slightly to get a map $\tilde{w}$ such that $\partial_1\tilde{w}^1 < 0$ on a small patch near the origin, and then minimize area among maps with the same boundary data as $\tilde{w}$ subject to the constraint that the first component is increasing in the $x_1$ direction. This gives rise to a free boundary problem (the free boundary region being that where the derivative in $x_1$ of the first component vanishes). We expect that solving it gives a map $w$ satisfying the conditions described before the statement of Theorem \ref{main}. This is the approach mentioned in Remark \ref{BigNHood}. However, it seems that this approach will require an involved analysis that we leave for another paper. We instead build $w$ by hand. More precisely, we use Cauchy-Kovalevskaya to build a map $w$ that satisfies the conditions described before the statement of Theorem \ref{main}, and in addition solves the Euler-Lagrange system associated to minimizing area subject to the gradient constraint $\partial_1w^1 \geq 0$.

Proving Theorem \ref{main} from here is still not straightforward. We need to use a new calibration argument, which leverages some convexity properties of the area integrand in higher codimension as well as the form of the Euler-Lagrange system for $w$, to show that the graph of $u$ has smaller area than all graphical competitors.

To conclude the introduction, we remark that there are some parallels between this construction and our construction of a non-$C^1$ solution to the special Lagrangian equation in a previous paper \cite{MS}. Namely, in that paper we first construct via Cauchy-Kovalevskaya a real-valued function $w$ on a domain in $\mathbb{R}^3$ that solves a degenerate Bellman equation (special Lagrangian equation outside some free boundary region, homogeneous Monge-Amp\`{e}re equation inside) and then take a Legendre transform (i.e. rotate its gradient graph in $\mathbb{R}^6$) to get the example. The gradient graph of that example also has a ``vertical" patch of dimension three that is non-minimal. However, from a variational perspective, that construction is fundamentally different from the one in this paper. We will illustrate this by showing that the gradient graph of the example from \cite{MS} does not minimize area among graphical competitors, and moreover not even among Lagrangian competitors.

The paper is organized as follows. In Section \ref{Prelims} we discuss some preliminaries. In particular, we review various ways of writing the minimal surface system, record a useful convexity property of the area integrand, and review a few facts from the theory of Cartesian currents. In Section \ref{PointSing} we construct the example with point singularity. In Section \ref{LargeSing} we construct the example with large singularity. Finally, in Section \ref{SLAG} we show that the special Lagrangian example from \cite{MS} is not area-minimizing among Lagrangian graphs.

\section*{Acknowledgments}
C. Mooney was supported by a Simons Fellowship, a Sloan Research Fellowship, NSF CAREER Grant DMS-2143668, and a UCI Chancellor's Fellowship. O. Savin was supported by NSF Grant DMS-2349794.

\section{Preliminaries}\label{Prelims}

\subsection{Minimal Surface System}
Let $u = (u^1,\,...,\,u^m)$ be a $C^2$ map from a domain $\Omega \subset \mathbb{R}^n$ to $\mathbb{R}^m$ and let $F$ be a smooth function on $\mathbb{R}^{m \times n}$. Assume that $u$ solves the system
$$\partial_i(F_{p^{\alpha}_i}(Du)) = f^{\alpha}, \quad \alpha = 1,\,...,\,m.$$
This is called the outer variation system, equivalent to 
$$E(u + \epsilon \varphi) = E(u) - \epsilon \int_{\Omega} f \cdot \varphi \,dx + o(\epsilon)$$
for all $\varphi \in C^1_0(\Omega)$, where 
$$E(v) := \int_{\Omega} F(Dv)\,dx.$$ 
Taking the dot product of the outer variation system with $\partial_ju$ we get
$$\partial_i(F(Du)\delta_{ij} - F_{p^{\alpha}_i}(Du)u^{\alpha}_j) = -\partial_j u \cdot f, \quad j = 1,\,...,\,n.$$
This is called the inner variation system, and it is equivalent to
$$E(u(\cdot + \epsilon V)) = E(u) - \epsilon\int_{\Omega} f \cdot (Du \cdot V)\,dx + o(\epsilon)$$ 
for any vector field $V \in C^{\infty}_0(\Omega)$. The derivation of the inner variation system is motivated by the fact that $u(x + \epsilon V) = u(x) + \epsilon Du(x) \cdot V(x) + o(\epsilon),$ 
i.e. that a domain variation by $V$ gives rise to an outer variation by $Du \cdot V$.

For the remainder of the paper we will let $F$ be the area integrand:
$$F(M) := \sqrt{\det(I + M^TM)}, \quad M \in \mathbb{R}^{m \times n}.$$
In terms of the principal values $\lambda_i$ of $M$, we have
\begin{equation}\label{PVals}
F(M) = \Pi_{i = 1}^n (1+\lambda_i^2)^{1/2}.
\end{equation}

For $u$ as above, we let 
$$g := I + Du^TDu.$$ 
Then 
$$F_{p^{\alpha}_i}(Du) = \sqrt{\det g}g^{ik}\partial_ku^{\alpha}$$ and 
$$\partial_ku \cdot \partial_ju = g_{kj} - \delta_{kj},$$ 
whence
$$F(Du)\delta_{ij} - F_{p^{\alpha}_i}(Du)u^{\alpha}_j = \sqrt{\det g}g^{ij}.$$
Thus, if
\begin{equation}\label{Outer}
\partial_i(\sqrt{\det g}g^{ij}\partial_ju^{\alpha}) = f^{\alpha},\, \alpha = 1,\,...,\,m
\end{equation}
then by the general considerations above we have
\begin{equation}\label{Inner}
\partial_i(\sqrt{\det g}g^{ij}) = -\partial_ju \cdot f, \, j = 1,\,...,\,n.
\end{equation}

This implication has a geometric interpretation. The mean curvature vector to the graph of $u$ is 
$$\Delta_g(x,\,u) := \frac{1}{\sqrt{\det g}}\partial_i(\sqrt{\det g}g^{ij}\partial_j(x,\,u)).$$ 
The outer variation system (\ref{Outer}) says that the vertical component of the mean curvature vector is $\frac{1}{\sqrt{\det g}}f$. The mean curvature vector is orthogonal to the graph of $u$, so the vertical component determines the horizontal one, given by the inner variation system (\ref{Inner}).

By combining the systems we see that if (\ref{Outer}) holds, then
\begin{equation}\label{MSS}
\sqrt{\det g}g^{ij}u^{\alpha}_{ij} = f^{\alpha} + (f \cdot \partial_ju)\partial_ju^{\alpha}, \quad \alpha = 1,\,...,\,m.
\end{equation}
The converse is in fact true. To see this, use again that the mean curvature vector is orthogonal to the graph of $u$:
\begin{align*}
0 &= \sqrt{\det g}\Delta_g(x,\,u) \cdot (e_k,\,\partial_ku) \\
&= [\partial_i(\sqrt{\det g}g^{ij})(e_j,\,\partial_ju) + (0,\,\sqrt{\det g}g^{ij}u_{ij})] \cdot (e_k,\,\partial_ku) \\
&\overset{(\ref{MSS})}{=} \partial_i(\sqrt{\det g}g^{ij})g_{jk} + f \cdot \partial_ku + (f \cdot \partial_ju)\partial_ju \cdot \partial_k u \\
&= [\partial_i(\sqrt{\det g}g^{ij}) + f \cdot \partial_ju]g_{jk}
\end{align*}
for all $k$. This immediately implies (\ref{Inner}), which combined with (\ref{MSS}) gives (\ref{Outer}). Below we will often find it convenient to use the non-divergence form (\ref{MSS}) of the system. When $f = 0$, we call it the minimal surface system.

\begin{rem}
In this remark we take $f = 0$. When $u$ is merely Lipschitz, one can still interpret the systems (\ref{Outer}) and (\ref{Inner}) in the sense of distributions. However, in that case, the implication (\ref{Outer}) $\Rightarrow$ (\ref{Inner}) is not clear. It was conjectured to be true by Lawson and Osserman (see \cite{LO}, Conjecture 2.1). In \cite{HMT} it was shown by Hirsch, the first-named author and Tione that Lipschitz solutions to (\ref{Outer}) are smooth when $n = 2$ and $m$ is arbitrary, resolving the conjecture for surfaces. The conjecture remains open when $n \geq 3,\, m \geq 2$. This problem is closely related to the ``wild" convex integration solutions to outer variation systems associated to strictly quasiconvex, resp. polyconvex integrands constructed by M\"{u}ller-\v{S}ver\'{a}k in \cite{MuS}, resp. Sz\'{e}kelyhidi in \cite{Sz}.
\end{rem}

\begin{rem}
If both (\ref{Outer}) and (\ref{Inner}) hold for a Lipschitz map $u$ with $f = 0$, then the graph of $u$ is a critical point of area with respect to all deformations of the ambient space $\mathbb{R}^{n+m}$ supported away from the boundary. Thus, the monotonicity formula holds. We note in particular that if a minimizer $u$ of the area functional $\int F(Du)$ existed in the class of Lipschitz maps, then it would solve both systems and the monotonicity formula would hold. Indeed, both outer and domain deformations would give rise to admissible competitors.
\end{rem}

\subsection{A Convexity Property}
We claim that there exists $\delta(n,\,m) > 0$ such that whenever $N,\, M \in \mathbb{R}^{m \times n}$ and $|N| < \delta$, 
\begin{equation}\label{ConvIneq}
F(M) \geq |F(N) + DF(N) \cdot (M-N)|.
\end{equation}

To see this we use two simple facts:
\begin{equation}\label{FExp}
F(M) = 1 + \frac{1}{2}|M|^2 + O(|M|^4) \text{ near } 0, \text{ and }
\end{equation}
\begin{equation}\label{FBound}
F(M) \geq (1+|M|^2)^{1/2}.
\end{equation}
These facts are easy to verify by looking at the expression (\ref{PVals}) for $F$ in terms of the principal values of $M$.

Using (\ref{FExp}) we see that $F$ is strictly convex in $B_c$ for some $c(n,\,m) > 0$. For $\delta < c$ to be chosen, let $\tilde{F}$ be the supremum of supporting linear functions to $F$ in $B_{\delta}$. Taking $\delta$ small we can guarantee that 
$$|D\tilde{F}| \leq 2\delta,$$ 
thus by taking $\delta$ smaller if necessary we have that
$$\tilde{F}(M) \leq 1 + 2\delta|M| < (1+|M|^2)^{1/2} \text{ outside of } B_c.$$
Thus, $\tilde{F}$ is convex, agrees with $F$ in $B_{\delta}$, and by the convexity of $F$ in $B_c$ and (\ref{FBound}) satisfies 
$$\tilde{F} \leq F$$
globally. 

For $N \in B_{\delta}$ we conclude that
$$F(M) \geq \tilde{F}(M) \geq \tilde{F}(N) + D\tilde{F}(N) \cdot (M-N) = F(N) + DF(N) \cdot (M-N).$$
To prove the other direction of (\ref{ConvIneq}), use that $\tilde{F} \geq 1$, the bound on $|D\tilde{F}|$, and (\ref{FBound}):
$$-\tilde{F}(N) - D\tilde{F}(N) \cdot (M-N) \leq -1 + 2\delta(|M| + \delta) \leq |M| < F(M).$$

\begin{rem}\label{StrictConvIneq}
By inspecting the proof, one can also see that the inequality (\ref{ConvIneq}) is strict unless $M = N$.
\end{rem}

\subsection{Cartesian Currents}\label{CartCurr}
An introduction to the theory of Cartesian currents can be found in Sections 2 and 3 of the paper \cite{GMS4} of Giaquinta-Modica-Sou\v{c}ek (see also Leon Simon's notes \cite{Si} for an introduction to the general theory of currents). We record a few facts here for the reader's convenience.

The space $D_n(U)$ of $n$-dimensional currents on a domain $U \subset \mathbb{R}^{n+m}$ is dual to the space $D^n(U)$ of smooth compactly supported $n$-forms on $U$. The boundary of a current $T \in D_n(U)$ is an $n-1$ dimensional current defined by enforcing Stokes' theorem: 
$$\partial T(\omega) = T(d\omega) \text{ for all } \omega \in D^{n-1}(U).$$ 
The mass of a current $T$ is defined to be 
$$\sup_{\{\omega \in D^n(U),\, |\omega| \leq 1\}} T(\omega),$$ 
and this is easily seen to be lower semicontinuous with respect to weak convergence (that is, $T_k(\omega) \rightarrow T(\omega)$ for all $\omega \in D^n(U)$).

A useful sub-class of currents consists of the (integer multiplicity) rectifiable currents, which are given by a countably $n$-rectifiable set $M$, a positive integer-valued and $\mathcal{H}^n$-integrable function $\theta$ representing multiplicity, and an $\mathcal{H}^n$-measurable orientation $\xi$ for the approximate tangent spaces to $M$. The rectifiable current 
$$T = (M,\,\theta,\,\xi)$$ 
is defined by integration of $n$-forms against $\theta\xi$ on $M$. An important compactness theorem of Federer-Fleming \cite{FF} is that any sequence of rectifiable currents $T_k$ such that both $T_k$ and $\partial T_k$ have locally uniformly bounded mass has a sub-sequence converging weakly to a rectifiable current. Combined with lower semicontinuity of mass, this gives easily the existence of mass-minimizers among rectifiable currents with a fixed e.g. smooth, compact boundary.

Finally, let $\Omega \subset \mathbb{R}^n$ be a smooth bounded domain, and let 
$$U = \Omega \times \mathbb{R}^m \subset \mathbb{R}^{n+m}.$$
The Cartesian currents are the subset of rectifiable currents given by weak limits of graphs of Lipschitz maps $u_k : \Omega \rightarrow \mathbb{R}^m$ with uniformly bounded mass and $L^1$ norm. It is shown in \cite{GMS4} that a Cartesian current $T = (M,\, \theta,\,\xi)$ satisfies that the set $M_+ \subset M$ where $M$ has a tangent plane with no vertical direction has a one-to-one projection to a set of full measure in $\Omega$ and that $\theta = 1$ a.e. on $M_+$. The current $(M_+,\,1,\,\xi)$ is in fact given by the graph of a BV map, and $T$ can be decomposed as a sum of this current and a purely vertical part. It is also shown in \cite{GMS4} that sequences of Cartesian currents with uniformly bounded mass and $L^1$ norm (defined as 
$$\sup_{|\phi|\leq 1, \phi \in C^{\infty}_0(U)} T(|y|\phi\,dx),$$ 
where $(x,\,y) \in \Omega \times \mathbb{R}^m$) have subsequences converging weakly (in the sense of currents) to a Cartesian current. The existence of mass-minimizers among Cartesian currents with a fixed boundary given e.g. by the graph of a smooth map on $\partial \Omega$ follows easily.

\section{Example with point singularity}\label{PointSing}
In this section we build a smooth minimal graph of codimension two in $\mathbb{R}^5$ that has a vertical tangent plane at a single point (the origin). The idea is to first construct, via Cauchy-Kovalevskaya, an analytic solution $(w^1,\,w^2)$ to the minimal surface system near the origin in $\mathbb{R}^3$ such that $\partial_1w^1$ has a non-degenerate local minimum at $0$. We then exchange the $x_1$ and $w^1$ axes. 

\begin{rem}\label{MP}
Note that any solution $w$ to the minimal surface system, in any dimension and codimension, satisfies (after differentiating (\ref{MSS}))
\begin{equation}\label{OneDiff}
g^{ij}w^1_{1ij} + \partial_1g^{ij}w^1_{ij} = 0.
\end{equation}
In dimension $n = 2$ or codimension $m = 1$, the second term in (\ref{OneDiff}) can be written in terms of derivatives of $w^1_1$. In particular, derivatives of $w$ obey the maximum principle in those cases. This shows that our example has the smallest possible dimension and codimension. 
\end{rem}

\begin{rem}
In the case $n = m = 3$ similar examples are known, due to Nadirashvili-Vladut (\cite{NV2}) and Wang-Yuan (\cite{WY2}), which moreover are special Lagrangian.
\end{rem}

Let $\epsilon > 0$ and let $(w^1,\,w^2)$ be the solution to the minimal surface system near the origin in $\mathbb{R}^3$ obtained by Cauchy-Kovalevskaya with the data
$$(w^1,\,w^2)|_{\{x_1 = 0\}} = (x_2x_3,\, \epsilon x_2), \quad (\partial_1w^1,\,\partial_1w^2)|_{\{x_1 = 0\}} = (x_2^2+x_3^2,\, A_{\epsilon}x_3).$$
Here 
$$A_{\epsilon} := \frac{5 + 4\epsilon^2}{\epsilon}.$$ 
We claim that
$$w^1 = x_2x_3 + x_1(x_1^2+x_2^2+x_3^2) + O(|x|^4), \quad w^2 = \epsilon x_2 + A_{\epsilon} x_1x_3 + O(|x|^3).$$
This can be checked by calculating the derivatives of $w$ at the origin, using the Cauchy data and the system (\ref{MSS}). Here are some details. First, we have
$$g^{-1}|_{\{x_1 = 0\}} =  
\left(\begin{array}{ccc}
1 & -\frac{\epsilon A_{\epsilon}}{1+\epsilon^2}x_3 & 0 \\
-\frac{\epsilon A_{\epsilon}}{1+\epsilon^2}x_3 & \frac{1}{1+\epsilon^2} & 0 \\
0 & 0 & 1
\end{array}\right) + O(|x|^2).
$$
Using the system (\ref{MSS}) and the Cauchy data for $w^1,\, w^2$ we conclude that 
$$w^1_{11}|_{\{x_1 = 0\}} = O(|x|^2), \quad w^2_{11}(0) = 0.$$
The expansion for $w^2$ follows immediately, as well as 
$$w^1 = x_2x_3 + O(|x|^3).$$ 
We also know all third derivatives of $w^1$ at $0$ except for $w^1_{111}(0)$. Among these, all vanish apart from $$w^1_{122}(0) = w^1_{133}(0) = 2.$$ 
To find $w^1_{111}(0)$ we use the once-differentiated equation (\ref{OneDiff}) and evaluate at the origin:
$$\left(w^1_{111}(0) + \frac{2}{1+\epsilon^2} + 2\right) + 2\partial_1g^{23}(0) = 0.$$
We have
\begin{align*}
&\partial_1g^{23}(0) = -g^{2i}g^{3j}\partial_1g_{ij}(0) = -\frac{1}{1+\epsilon^2}\partial_1g_{23}(0) \\
&= -\frac{1}{1+\epsilon^2} (w^{\gamma}_{12}w^{\gamma}_3 + w^{\gamma}_2w^{\gamma}_{13})(0) = -\frac{1}{1+\epsilon^2} w^2_2w^2_{13}(0) = -\frac{\epsilon A_{\epsilon}}{1+\epsilon^2},
\end{align*}
and using the choice of $A_{\epsilon}$ we conclude that 
$$w^1_{111}(0) = 6,$$ 
as desired.

We now complete the construction. Note that 
$$w^1_1 = 3x_1^2 + x_2^2 + x_3^2 + O(|x|^3),$$
so in a neighborhood of the origin, $w^1_1$ is nonnegative and vanishes only at the origin. As a result, the map
$$y(x) = (w^1(x),\,x_2,\,x_3)$$
is injective in a neighborhood of the origin, and an analytic local diffeomorphism away from the origin. The graph of $(w^1,\,w^2)$ over a small ball $B$ can thus be represented as the graph of another map $(u^1,\,u^2)$ over $y(B)$ in the following way:
$$\{(x_1,\,x_2,\,x_3,\,w^1(x),\,w^2(x))\} = \{(u^1(y),\,y_2,\,y_3,\,y_1,\,u^2(y))\},$$
where 
$$u^1(y(x)) := x_1, \quad u^2(y(x)) := w^2(x).$$ 
The graph of $(u^1,\,u^2)$ is a rotation of an analytic minimal graph, so $(u^1,\,u^2)$ solves the minimal surface system where it is smooth (away from the origin). The relation
$$u^1(x_1^3 + O(x_1^4),\,0,\,0) = x_1$$
shows that $u$ is no better than $C^{1/3}$, and that $u^1_1$ tends to infinity near the origin, where the graph of $u$ has a tangent $3$-plane with one vertical direction.

To conclude we note that if we take $\epsilon$ and $B$ sufficiently small, then $|Dw|$ is very small and the graph of $w$ over $B$ is area-minimizing in $\mathbb{R}^5$ (see e.g. the work of Lawlor-Morgan \cite{LM} for a proof using the slicing method, which works even for unoriented competitors; see also work of Federer \cite{F} for a proof for oriented competitors, as well as Section \ref{Calibration} below for an elementary proof for oriented competitors using the calibration method). In particular, the graph of $u$ over $y(B)$ has smaller area than all graphical (in $y$) competitors with the same boundary. Thus, minimizing area among Lipschitz graphs can lead to an interior vertical tangent plane in dimension $n \geq 3$ and codimension $m \geq 2$.

\section{Example with large singularity}\label{LargeSing}
The example in the previous section had a vertical tangent plane at a single point. In this section we modify the construction to get an example that is vertical on a whole three-dimensional patch. The approach is to first minimize area subject to the convex gradient constraint $w^1_1 \geq 0$. This leads to an interesting free boundary problem. We then swap the $x_1$ and $w^1$ axes as before.

\subsection{A Free Boundary Problem}\label{FBP}
Suppose that $w = (w^1,\,w^2,\,...,\,w^m)$ is a $C^{1,1}$ map from a bounded domain $\Omega \subset \mathbb{R}^n$ to $\mathbb{R}^m$, such that $w^1_1 \geq 0$. Assume further that there is a smooth function $H : \Omega \rightarrow \mathbb{R}$ such that
\begin{equation}\label{EL0}
\{w^1_1 = 0\}  = \{H \geq 0\} \subset \subset \Omega
\end{equation}
and
\begin{equation}\label{ELSyst}
\partial_i(\sqrt{\det g}g^{ij}\partial_jw^{\alpha}) = \begin{cases}
\partial_1H\chi_{\{H > 0\}}, \quad \alpha = 1, \\
0, \quad \alpha > 1.
\end{cases}
\end{equation}
We claim that $w$ solves the Euler-Lagrange system associated to the area functional subject to the gradient constraint $w^1_1 \geq 0$. Indeed, let $\varphi \in C^{\infty}_0(\Omega)$ satisfy that $\partial_1\varphi^1 \geq 0$ in $\{H \geq 0\}$. We need to verify that
$$-\int_{\Omega} \partial_i(F_{p^{\alpha}_i}(Dw))\varphi^{\alpha} \geq 0.$$
Using (\ref{ELSyst}) this reduces to
$$-\int_{\{H > 0\}} \partial_1H\varphi^1 \geq 0,$$
which is easily verified by integrating by parts.

\begin{rem}\label{MP2}
Such a map cannot exist in the cases $n = 2$ or $m = 1$. Indeed, suppose one did. The form (\ref{MSS}) of the system gives
$$\sqrt{\det g}g^{ij}w^1_{ij} = \partial_1H\chi_{\{H > 0\}}(1+|\nabla w^1|^2).$$ 

When $n = 2$ this gives in $\{H > 0\}$ (using that $w^1_1 \equiv 0$ there) that 
$$\sqrt{\det g}g^{22}w^1_{22} = \partial_1H(1+(w^1_2)^2), \quad \partial_1(\sqrt{\det g}g^{22})w^1_{22} = \partial_1^2H(1+(w^1_2)^2).$$ 
These equations imply that $\partial_1H/(\sqrt{\det g}g^{22})$ is locally constant on lines in the $x_1$ direction in $\{H > 0\}$. Since $\sqrt{\det g}g^{22}$ is positive, we see that $\{H > 0\}$ cannot be compactly contained in $\Omega$. 

When $m = 1$, the left hand side of the equation is constant in the $x_1$ direction in $\{H > 0\}$. Differentiating the equation in the $x_1$ direction gives 
$$\partial_1^2H = 0 \text{ in }\{H > 0\},$$ 
leading to the same contradiction.
\end{rem}

We will now construct such a map in the minimal possible dimension and codimension $n = 3$ and $m = 2$. Let $\epsilon > 0$ (we'll fix it later). Let 
$$v^1 = x_2x_3,$$ 
and let $v^2$ be the solution to
$$\partial_i(\sqrt{\det g_v}g_v^{ij}\partial_jv^2) = 0$$
in a neighborhood of the origin with the Cauchy data 
$$v^2|_{\{x_1 = 0\}} = \epsilon x_2, \quad \partial_1v^2|_{\{x_1 = 0\}} = x_3.$$ 
Here $g_v = I + Dv^TDv$. Let 
$$f := \partial_i(\sqrt{\det g_v}g_v^{ij}\partial_jv^1).$$
The form (\ref{MSS}) of the system for $v$ gives an expression for $f$:
$$f = \frac{2\sqrt{\det g_v}g_v^{23}}{1+x_2^2+x_3^2}.$$
We calculate using the Cauchy data that 
$$g_v^{23}|_{\{x_1 = 0\}} = O(|x|^2).$$ 
We also calculate in a similar way to the previous section that
$$\partial_1g_v^{23}(0) = -g_v^{2i}g_v^{3j}\partial_1[(g_v)_{ij}](0) = -\frac{1}{1+\epsilon^2}\partial_1[(g_v)_{23}](0) = -\frac{\epsilon}{1+\epsilon^2}.$$
Thus, 
$$g_v^{23} = -\frac{\epsilon}{1+\epsilon^2}x_1 + O(|x|^2).$$ 
Combining this with the easy calculation 
$$\det g_v = 1+\epsilon^2 + O(|x|)$$ 
we conclude that
$$f = -\frac{2\epsilon}{\sqrt{1+\epsilon^2}} x_1 + O(|x|^2) = -C_{\epsilon}x_1 + O(|x|^2), \quad C_{\epsilon} := \frac{2\epsilon}{\sqrt{1+\epsilon^2}}.$$
Writing 
$$f = -C_{\epsilon}x_1 + h(x)$$ 
and taking 
$$H = \delta - C_{\epsilon}|x|^2/2 + \int_0^{x_1} h(t,\,x_2,\,x_3)\,dt = \delta - C_{\epsilon}|x|^2/2 + O(|x|^3),$$
we get 
$$f = \partial_1 H.$$ 
Moreover, for $\delta > 0$ small, the connected component $\Omega_0$ of $\{H > 0\}$ containing the origin is analytic and uniformly convex.

Now, let $w = v$ in $\overline{\Omega_0}$, let $\tilde{v}$ be the solution to the minimal surface system in a neighborhood of $\partial \Omega_0$ with the same values and derivatives as $v$ on $\partial \Omega_0$ (obtained by Cauchy-Kovalevskaya), and let $w = \tilde{v}$ outside of $\overline{\Omega_0}$. We will show that 
$$w^1_1 > 0 \text{ outside of } \overline{\Omega_0}.$$
In particular, $w$ satisfies (\ref{EL0}) and (\ref{ELSyst}).

Let $\Gamma \subset \partial \Omega_0$ be the analytic curve $\{\partial_1H = 0\} \cap \partial \Omega_0$, and let $\nu$ be the exterior unit normal to $\partial \Omega_0$. It suffices to show that
$$\partial_{\nu} \tilde{v}^1_1 > 0 \text{ on } \partial \Omega_0 \backslash \Gamma, \quad \partial_{\nu}^2\tilde{v}^1_1 > 0 \text{ on } \Gamma.$$
To that end we use the equations
$$\sqrt{\det g_v}g_v^{ij}v^{1}_{ij} = (1+x_2^2+x_3^2)\partial_1H, \quad \sqrt{\det g_{\tilde{v}}}g_{\tilde{v}}^{ij}\tilde{v}^{1}_{ij} = 0.$$
Here $g_v$, resp. $g_{\tilde{v}}$ denote the metrics $I + Dv^TDv$, resp. $I + D\tilde{v}^TD\tilde{v}$, and the first equation is obtained by looking at the form (\ref{MSS}) of the system for $v$. The metrics match on $\partial \Omega_0$, where we will drop the subscripts. Pick a point $p$ on $\partial \Omega_0$, and choose coordinates where $\nu$ at that point is a coordinate direction. All second derivatives of $v^1,\,\tilde{v}^1$ in these coordinates match except for the pure normal second derivatives. Subtracting the equations gives
$$v^1_{\nu\nu} - \tilde{v}^1_{\nu\nu} = \left(\frac{1+x_2^2+x_3^2}{\sqrt{\det g}g^{\nu\nu}}\right)\partial_1H$$
at $p$. If $p$ is not in $\Gamma$, then $\nu = A e_1 + B\tau$, where $A \neq 0$ has the same sign as $-\partial_1H$ and $\tau$ is tangent to $\partial \Omega_0$ at $p$. Using that $v^1_1 \equiv 0$ and that $v^1_{\nu\tau} = \tilde{v}^1_{\nu\tau}$ at $p$ we conclude from the above identity that
$$\tilde{v}^1_{1\nu} = -\left(\frac{1+x_2^2+x_3^2}{\sqrt{\det g}g^{\nu\nu}}\right)\frac{\partial_1H}{A} > 0.$$

Now suppose that $p \in \Gamma$. Since $\partial_1H = 0$ on $\Gamma$, the equations imply that 
$$D^2v^1 = D^2\tilde{v}^1 \text{ on } \Gamma.$$
In addition, since $e_1$ is tangent to $\partial \Omega_0$ at $p$, the derivative in the $e_1$ direction of the metrics agree at $p$. Differentiating the equations in the $e_1$ direction, subtracting, and evaluating at $p$, we thus get
$$\sqrt{\det g}g^{ij}(\tilde{v}^1_{1ij} - v^1_{1ij})(p) = -(1+x_2^2+x_3^2)\partial_1^2H(p) > 0.$$
Choose coordinates at $p$ where one coordinate direction is $\nu$, another is $\tau$ tangent to $\Gamma$, and the third is $e = \alpha e_1 + \beta \tau$ for some $\alpha > 0$. Since the Hessians of $v$ and $\tilde{v}$ match on $\Gamma$, all third derivatives involving a differentiation in $\tau$ agree at $p$, giving at $p$ that
$$g^{ee}(\tilde{v}^1_{1ee} - v^1_{1ee}) + 2\alpha g^{e\nu}(\tilde{v}^1_{\nu11}-v^1_{\nu11}) + g^{\nu\nu}(\tilde{v}^1_{1\nu\nu}-v^1_{1\nu\nu}) > 0.$$
Since $e$ and $e_1$ are tangent to $\partial \Omega_0$ at $p$ and the gradients of $\tilde{v}^1$ and $v^1$ agree on $\partial \Omega_0$, the first and second terms can be rewritten as expressions involving the curvature of $\partial \Omega_0$ times differences between second derivatives of $\tilde{v}^1$ and $v^1$ at $p$. Using again that $D^2v^1 = D^2\tilde{v}^1$ at $p$ we see that the first two terms vanish at $p$. Since $v^1_1 \equiv 0$ we conclude that 
$$\tilde{v}^1_{1\nu\nu}(p) > 0,$$ 
as desired.

\begin{rem}\label{PP}
By taking $\epsilon$, then $\delta$ small, we can make $|Dw|$ as small as we like in $\Omega$ (here $\Omega$ is a small convex neighborhood of $\overline{\Omega_0}$). We can also ensure that there is an analytic solution $\tilde{w}$ to the minimal surface system in $\Omega$ with the same boundary data as $w$, with $|D\tilde{w}|$ as small as we like. This can be accomplished e.g. via the implicit function theorem in Banach spaces, using that the minimal surface system linearizes to the Laplace operator in each component at $0$, or by appealing to the non-perturbative existence result proved by M.-T. Wang in \cite{W}. The graph of $\tilde{w}$ is area-minimizing for the same reasons cited at the end of the previous section, and $\partial_1\tilde{w}^1 < 0$ somewhere, in view of the calibration arguments later in this section.
\end{rem}

\subsection{Large Singularity}\label{LargeSingSub}
We again let 
$$y(x) = (w^1(x),\,x_2,\,x_3).$$
Then $y\left(\overline{\Omega_0}\right)$ is an analytic surface $\Sigma$ with boundary that resembles a potato chip, and lies in $\{(x_2x_3,\,x_2,\,x_3)\}$. Moreover, since $\partial_1w^1 > 0$ outside $\overline{\Omega_0}$, we see that $y$ maps $\Omega \backslash \overline{\Omega_0}$ diffeomorphically onto $U \backslash \Sigma$, where 
$$U := y(\Omega)$$ 
is a neighborhood of $\Sigma$ in $\mathbb{R}^3$. We can again represent the graph of $w$ over $\Omega \backslash \overline{\Omega_0}$ by the graph of an analytic map $u$ on $U \backslash \Sigma$ defined by
$$u^1(y(x))= x_1, \quad u^2(y(x)) = w^2(x).$$
Here and below we in fact view the graph of $u$ as a rotation of the graph of $w$. That is, we view the graph of $u$ has having a three-dimensional vertical patch that fills the ``hole" in the graph of $u$ over $U \backslash \Sigma$. Viewed this way, each point in the interior of $\Sigma$ is the projection to the $y$ subspace of an analytic one-dimensional curve that lies in the $u^1-u^2$ subspace.

\subsection{Calibration}\label{Calibration}
Here we aim to show that the graph of $u$ (interpreted as above, i.e. as a rotation of the graph of $w$) has smaller area than the graph of any Lipschitz map from $U$ to $\mathbb{R}^2$ that agrees with $u$ on $\partial U$. This illustrates that, when we minimize area while imposing a graphicality condition in higher codimension, large vertical and non-minimal portions can develop {\it away from the boundary}, in analogy with what happens at the boundary in the codimension one case.

We begin with a general observation. Let $(x_1,\,...,\,x_n,\,z_1,\,...,\,z_m)$ be coordinates on $\mathbb{R}^{n+m}$, and let $w = (w^1,\,...,\,w^m)$ be a $C^1$ map from a bounded domain $\Omega \subset \mathbb{R}^n$ into $\mathbb{R}^m$. Let $F$ be the area integrand. We define the $n$-form
$$\omega := (F(Dw) - F_{p^{\alpha}_i}(Dw)w^{\alpha}_i)dx_1 \wedge ... \wedge dx_n + F_{p^{\alpha}_i}(Dw) \omega^{\alpha}_i$$
on $\Omega \times \mathbb{R}^m,$ where 
$$\omega^{\alpha}_i = dx_1 \wedge ... \wedge dx_{i-1} \wedge dz_{\alpha} \wedge dx_{i+1} \wedge ... dx^n.$$
We first claim that 
$$|\omega(T)| \leq 1$$ 
for any unit volume $n$-plane $T$, provided $|Dw| < \delta(n,\,m)$ small. Indeed, to check this on $n$-planes that project non-singularly to $x$, it suffices to do so for 
$$T = F(M)^{-1}[(e_{\sigma(1)} + v_{\sigma(1)}) \wedge ... \wedge (e_{\sigma(n)} + v_{\sigma(n)})],$$
where $\sigma$ is any permutation of $\{1,\,...,\,n\}$, $e_i$ are the unit coordinate vectors in the $x$ subspace, $v_i$ are vectors in the $z$ subspace, and $M \in \mathbb{R}^{m \times n}$ is the matrix whose $i^{th}$ column is $v_i$. Direct calculation gives
$$\omega(T) = \text{sgn}(\sigma) F(M)^{-1}[F(Dw) + DF(Dw) \cdot (M-Dw)],$$
and $|\omega(T)| \leq 1$ by virtue of the convexity inequality (\ref{ConvIneq}). The inequality for arbitrary unit $n$-planes follows by approximation (or by noting that the first term in $\omega$ gives no contribution for $n$-planes with a vertical direction, and $|F_{p^{\alpha}_i}(Dw)| \leq C(n,\,m)\delta$). The above calculation and Remark \ref{StrictConvIneq} also show that $\omega(T) = 1$ at $(x,\,z)$ only when $T$ is the tangent $n$-plane to the graph of $w$ at $x$.

We now observe that 
$$d\omega = -\partial_i(F_{p^{\alpha}_i}(Dw))dz_{\alpha} \wedge dx_1 \wedge ... \wedge dx_n.$$
In particular, $\omega$ is closed if $w$ solves the minimal surface system. Hence, provided $|Dw|$ is small in $\Omega$, the form $\omega$ is a calibration in the cylinder over $\Omega$ (so the graph of $w$ minimizes area among competing currents with support in $\Omega \times \mathbb{R}^m$). If $\Omega$ is convex, then the nearest point projection to $\Omega \times \mathbb{R}^m$ is area-decreasing, and we conclude that the graph of $w$ is area-minimizing in $\mathbb{R}^{n+m}$.

Finally, we specialize to maps $w$ solving the free boundary problem (\ref{EL0}), (\ref{ELSyst}), and satisfying in addition that $|Dw|$ is small (as the specific example we constructed in the case $n = 3,\,m = 2$ does, see Remark \ref{PP}). Let 
$$y(x) = (w^1(x),\,x_2,\,...,\,x_n),$$
and define a map $u = (u^1,\,...,\,u^m)$ on $y(\Omega)$ by 
$$u^1(y(x)) = x_1, \quad u^{\alpha}(y(x)) = w^{\alpha}(x) \text{ for } \alpha > 1.$$ 
Here we are abusing notation; $y$ is only a diffeomorphism of $\{H < 0\}$, so $u$ is only a smooth map on $y(\{H < 0\})$. As usual, we view the graph of $u$ as a rotation of the graph of $w$:
$$\{(x_1,\,...,\,x_n,\,w^1(x),\,...,\,w^m(x))\} = \{(u^1(y),\,y_2,\,...,\,y_n,\,y_1,\,u^2(y),\,...,\,u^m(y))\},$$
and we note $u$ is a smooth map near $\partial (y(\Omega))$. We have by (\ref{ELSyst}) that
$$d\omega = -\partial_i(F_{p^{\alpha}_i}(Dw))dz_{\alpha} \wedge dx_1 \wedge ... \wedge dx_n = d \tilde{\omega},$$
where 
$$\tilde{\omega} = H(x)\chi_{\{H > 0\}} dz_1 \wedge dx_2 \wedge ... \wedge dx_n.$$
Since $w^1_1 = 0$ in $\{H > 0\}$, $\tilde{\omega}$ vanishes on the tangent planes to the graph of $w$. Let $\Gamma_{\psi}$ be the graph of a Lipschitz map $\psi$ of $(z_1,\,x_2,\,...,\,x_n)$ with the same values as $u$ on $\partial (y(\Omega))$, and assume first that $\Gamma_{\psi}$ is contained in the cylinder $\Omega \times \mathbb{R}^m$. Let $\Gamma_u$ be the graph of $u$. The above considerations and Stokes' theorem give
\begin{align*}
\text{vol}(\Gamma_{\psi}) - \text{vol}(\Gamma_u) &\geq \int_{\Gamma_{\psi}} \omega - \int_{\Gamma_u} \omega \quad (|\omega| \leq 1,\, = 1 \text{ on } \Gamma_u) \\
&= \int_{\Gamma_{\psi}} \tilde{\omega} - \int_{\Gamma_u} \tilde{\omega} \quad \text{ (Stokes' theorem) } \\
&= \int_{\Gamma_{\psi} \cap (\{H > 0\} \times \mathbb{R}^m)} H(x) dz_1 \wedge dx_2 \wedge ... \wedge dx_n \\
& \geq 0,
\end{align*}
where in the last inequality we use that $\Gamma_{\psi}$ carries the orientation of $\mathbb{R}^n$ with the coordinates $(z_1,\,x_2,\,...,\,x_n)$. In particular, the graph of $u$ has smaller area than graphical competitors that also lie in $\Omega \times \mathbb{R}^m$. If $\Omega$ is in addition convex (as it is for our specific example in the case $n = 3,\, m = 2$), we can use that the nearest point projection to $\Omega \times \mathbb{R}^m$ is area-decreasing and argue similarly (replacing $\Gamma_{\psi}$ by its projection to $\Omega \times \mathbb{R}^m$) to prove that the graph of $u$ has smaller area than any Lipschitz graph over $y(\Omega)$ with the same boundary.

\begin{rem}\label{CC}
The same argument shows that the graph of $u$ in fact minimizes mass among Cartesian currents with the same boundary, using that such currents carry the orientation of $\mathbb{R}^n = \{(z_1,\,x_2,\,...,\,x_n)\}$ (see e.g. \cite{GMS4}, Theorem 2 on pg. 108).
\end{rem}

\section{Special Lagrangian Example}\label{SLAG}
In \cite{MS} we constructed a viscosity solution to the special Lagrangian equation in a domain in $\mathbb{R}^3$ whose gradient graph contains a $3$-dimensional vertical portion that is not minimal. It is also constructed by solving a free boundary problem. Here we briefly recall the construction, and we show that the gradient graph of that example does not minimize area among graphical competitors (not even Lagrangian ones), showing that the example from \cite{MS} is different from the one in this paper from a variational perspective. See Remark \ref{Variational} below for a deeper discussion.

For a symmetric $n \times n$ matrix $M$ with eigenvalues $\{\lambda_i\}_{i = 1}^n$ we let 
$$G(M) := \sum_{i = 1}^n \tan^{-1}(\lambda_i).$$
For $\lambda > 0$ to be fixed later let
$$\Phi(x) := \frac{\lambda x_1^2}{1+x_3} + \frac{\lambda x_2^2}{1-x_3}.$$
Each term in $\Phi$ is a translation of a one-homogenous function of two variables, so $D^2\Phi$ has rank two, and $\det D^2\Phi = 0$.
We proved in \cite{MS} that the Lagrangian angle 
$$\Theta := G(D^2\Phi)$$ 
satisfies
\begin{equation}\label{ThetaExp}
\Theta - \Theta(0) = 2\lambda(x_1^2 + x_2^2 + 2x_3^2) + O(\lambda^2)O(|x|^2).
\end{equation}
We choose $\lambda > 0$ small so that $\Theta$ is uniformly convex near $0$. Then for $\epsilon > 0$ small the connected component $\Omega_0$ of $\{\Theta < \Theta(0) + \epsilon^2\}$ containing $0$ is an analytic, uniformly convex set. We defined $w$ to be $\Phi$ in $\overline{\Omega_0}$, and for a small neighborhood $\Omega$ of $\overline{\Omega_0}$ we define $w$ in $\Omega \backslash \overline{\Omega_0}$ to be the solution to 
$$G(D^2w) = \Theta(0) + \epsilon^2$$ 
with the same first-order expansion as $\Phi$ on $\partial \Omega_0$ (obtained by Cauchy-Kovalevskaya). We proved that 
$$\det D^2w < 0 \text{ outside of } \overline{\Omega_0}$$ 
(so $w$ solves the degenerate Bellman equation 
$$\max\{G(D^2w) - \Theta(0) - \epsilon^2,\, \det D^2w\} = 0),$$
and combining this with the information that
\begin{equation}\label{Hess}
D^2w = 2\lambda(I - e_3 \otimes e_3) + O(\epsilon)
\end{equation}
in $\Omega$ we proved that $\nabla w$ is a diffeomorphism of $\Omega \backslash \overline{\Omega_0}$ (and similarly to above, maps $\Omega_0$ to a potato chip-like analytic surface compactly contained in $\nabla w(\Omega)$). 

We claim that provided $\lambda$, then $\epsilon$ were chosen small enough, we can find a function $\varphi \in C^{\infty}_0(\Omega)$ such that, for all $t > 0$ small, $\nabla(w + t\varphi)$ is a diffeomorphism of $\Omega$ whose graph has smaller volume than that of $\nabla w$. This confirms the non-minimality stated at the beginning of this section, after taking Legendre transforms. More precisely, we view the gradient graph of the Legendre transform $w^*$ of $w$ as being a rotation of the gradient graph of $w$, with a three-dimensional vertical portion (a rotation of the graph of $\nabla w$ over $\Omega_0$) over the surface $\nabla w(\Omega_0)$. We have for $t > 0$ small that $(w + t\varphi)^*$ is well-defined, and has $C^{1,\,1}$ gradient graph competing with that of $w^*$ but with smaller volume.

To prove the claim we let $\rho \in C^{\infty}_0(\Omega)$ be a cutoff function that is $1$ in a neighborhood $N$ of $\overline{\Omega_0}$, and we take 
$$\varphi = (x_1^2+x_2^2-x_3^2)\rho.$$
Using (\ref{Hess}) we have
\begin{align*}
\det D^2(w + t\varphi) &= \det D^2w + (4\lambda^2\varphi_{33} + O(\epsilon)O(|D^2\varphi|))t \\
&+ O(\lambda)O(|D^2\varphi|^2)t^2 + O(|D^2\varphi|^3)t^3.
\end{align*}
We conclude that 
$$\det D^2(w + t\varphi) < 0$$
for $t > 0$ small, using that $\varphi_{33} = -2$ in $N$, that $|\varphi_{ij}| \leq 2$ in $N$, and that 
$$\det D^2w \leq -\kappa < 0$$ 
outside $N$. Similar arguments as in Lemma 2.4 of \cite{MS}, using (\ref{Hess}), imply that $\nabla (w + t\varphi)$ is a diffeomorphism of $\Omega$ for $t > 0$ small.

The first variation of the area of the graph of $\nabla w$ in direction $\nabla \varphi$ is given by 
$$\int \sqrt{\det g}g^{ij}w_{kj}\varphi_{ki}.$$
Here 
$$g = I + (D^2w)^2.$$ 
Since the graph of $\nabla w$ is minimal outside $\Omega_0$, we just need to show that
$$\int_{\Omega_0} \partial_i(\sqrt{\det g}g^{ij}w_{kj})\varphi_k > 0.$$
The expression $\partial_i(\sqrt{\det g}g^{ij}\partial_j\nabla w)$ is $\sqrt{\det g}$ times the vertical component of the mean curvature vector to the graph $\Sigma$ of $\nabla w$. The mean curvature vector of $\Sigma$ is given by the appealing formula 
$$J\nabla_{\Sigma}\Theta,$$
where
$$J(x,\,y) = (-y,\,x),\, x,\, y \in \mathbb{R}^3$$ 
and $\nabla_{\Sigma}$ is the gradient operator on $\Sigma$, as shown by Harvey-Lawson in \cite{HL1}. Since the tangent plane to $\Sigma$ has slope $D^2w = O(\lambda)$, the expansion (\ref{ThetaExp}) gives
$$J\nabla_{\Sigma}\Theta = (0,\,0,\,0,\,4\lambda x_1,\,4\lambda x_2,\,8\lambda x_3) + O(\lambda^2)O(|x|).$$
Since $\det g = 1 + O(\lambda^2)$ we conclude that 
\begin{equation}\label{MeanC}
\partial_i(\sqrt{\det g}g^{ij}\partial_j\nabla w) = (4\lambda x_1,\,4\lambda x_2,\,8\lambda x_3) + O(\lambda^2)O(|x|).
\end{equation}
Using (\ref{MeanC}), recalling the facts that $\varphi = x_1^2+x_2^2-x_3^2$ in $\Omega_0$ and
$$\Omega_0 = \{2\lambda(x_1^2 + x_2^2 + 2x_3^2) + O(\lambda^2)O(|x|^2) < \epsilon^2\},$$
and making the change of variable 
$$(x_1,\,x_2,\,x_3) = (2\lambda)^{-1/2}\epsilon(y_1,\,y_2,\,y_3/\sqrt{2}),$$ 
we conclude that
\begin{align*}
\int_{\Omega_0} \partial_i(\sqrt{\det g}g^{ij}w_{kj})\varphi_k\,dx &= 8\lambda \int_{\Omega_0} (x_1^2+x_2^2-2x_3^2 + O(\lambda)O(|x|^2))\,dx \\
&=\epsilon^5\lambda^{-3/2}\int_{\{|y|^2 + O(\lambda) < 1\}} (y_1^2+y_2^2-y_3^2 + O(\lambda))\,dy \\
&> 0,
\end{align*}
as desired. 

\begin{rem}\label{Variational}
A closer analogue to what we did in this paper would be to minimize area among graphs of maps $v: \Omega \rightarrow \mathbb{R}^3$ with the same boundary data as $\nabla w$, subject to the constraint $\det Dv \leq 0$. This section shows that the minimizer is most likely not Lagrangian. On the other hand, minimizing area among graphs of Lagrangian maps satisfying this determinant constraint leads to a free boundary problem for the Hamiltonian stationary equation, a fourth-order PDE which says that the Lagrangian angle is harmonic on the graph away from the free boundary region. In particular, the Lagrangian angle need not be constant, which is consistent with our finding here.
\end{rem}



\end{document}